\documentclass{amsart}
\usepackage[headings]{fullpage}
\usepackage{amssymb,epic,eepic,epsfig,amsbsy,amsmath}

\theoremstyle{plain}

\newtheorem{theorem}{Theorem}
\newtheorem{proposition}{Proposition}[section]
\newtheorem{lemma}[proposition]{Lemma}
\newtheorem{corollary}[proposition]{Corollary}

\theoremstyle{definition}

\theoremstyle{remark}

\newtheorem{remark}[proposition]{Remark}

\def\printname#1{
    \if\draft y
        \smash{\makebox[0pt]{\hspace{-0.5in}
            \raisebox{8pt}{\tt\tiny #1}}}
    \fi
}

\newlength{\standardunitlength}
\catcode`\@=11
\long\def\@makecaption#1#2{%
     \vskip 10pt

\setbox\@tempboxa\hbox{
       \small\sf{\bfcaptionfont #1. }\ignorespaces #2}%
     \ifdim \wd\@tempboxa >\captionwidth {%
         \rightskip=\@captionmargin\leftskip=\@captionmargin
         \unhbox\@tempboxa\par}%
       \else
         \hbox to\hsize{\hfil\box\@tempboxa\hfil}%
     \fi}
\font\bfcaptionfont=cmssbx10 scaled \magstephalf
\newdimen\@captionmargin\@captionmargin=2\parindent
\newdimen\captionwidth\captionwidth=\hsize
\catcode`\@=12

\newcommand{\tr}{\operatorname{tr}}

\newcommand{\id}{\operatorname{id}}
\newcommand{\Vol}{\operatorname{Vol}}

\newcommand{\qbinom}[2]{\text{$\left[\begin{array}{c}#1\\ #2\end{array}
\right]$}}

\def\lbl#1{\label{#1}\printname{#1}}

\def\BZ{\mathbb Z}

\def\BR{\mathbb R}
\def\BC{\mathbb C}

\def\bb{\mathfrak b}

\def\A{\mathcal A}

\def\U{\mathcal U}
\def\D{\Delta}

\def\RR{{\mathcal P}}

\def\J{\mathcal J}

\def\F{\mathcal F}

\def\P{\mathcal P}
\def\R{\mathcal R}

\def\I{\mathcal I}

\def\g{\frak g}

\def\g{\mathfrak g}
\def\fg{\mathfrak g}

\def\ve{\varepsilon}

\newcommand{\inv}{\operatorname{inv}}

\def\Det{\widetilde{\det}}

\begin{document}


\title[Colored Jones Polynomial ]{
On the Colored Jones Polynomial and the Kashaev invariant}

\author{Vu Huynh \& Thang T. Q. L\^e}
\address{Department of Mathematics\\
SUNY Buffalo\\
Buffalo, NY 14260, USA}
\address{School of Mathematics \\
         Georgia Institute of Technology \\
         Atlanta, GA 30332-0160, USA}

\email{letu@math.gatech.edu}

\thanks{The second author was supported in part by National Science Foundation. \\
\newline
1991 {\em Mathematics Classification.} Primary 57M25.
\newline
}

\date{\today \hspace{0.5cm} First edition: January 31, 2005.}


\begin{abstract}
We express the colored Jones polynomial as the inverse of the
quantum determinant of a matrix with entries in the $q$-Weyl
algebra of $q$-operators, evaluated at the trivial function (plus
simple substitutions). The Kashaev invariant is proved to be equal
to  another special evaluation of the determinant. We also discuss
the similarity between our determinant formula of the Kashaev
invariant and the determinant formula of the hyperbolic volume of
knot complements, hoping it would lead to a proof of the volume
conjecture.
\end{abstract}

\maketitle


\addtocounter{section}{-1}


\section{Introduction}

For a knot $K$ in  $\BR^3$, the colored Jones polynomial $J'_K(N)$
is a Laurent polynomial, $J'_K(N) \in \R := \BZ[q^{\pm1}]$, see
\cite{Jones,MMe}. Here $N$ is a positive integer standing for the
$N$-dimensional prime $sl_2$-module. We use the unframed version
and the normalization in which $J'_K(N) = 1$ when $K$ is the
unknot. The colored Jones polynomial $J'_K(N)$ is defined using
the $R$-matrix of the quantized enveloping algebra of $sl_2(\BC)$.

Here we present the colored Jones polynomial as the inverse of the
quantum determinant of an almost quantum matrix whose entries are
in the  $q$-Weyl algebra of $q$-operators acting on the polynomial
rings, evaluated at the constant function 1. The proof is based on
the quantum MacMahon Master theorem proved in \cite{GLZ}.
Actually, it was an attempt to get a determinant formula for the
colored Jones polynomial that led the second author to the
conjecture that eventually became the quantum MacMahon's Master
theorem in \cite{GLZ}.

We will then give an application to the case of the Kashaev
invariant $\langle K\rangle_N := J'_K(N)|_{q=\exp 2\pi i/N}$. We
show that a special evaluation of the determinant will give the
Kashaev invariant. Our interpretation of the Kashaev invariant
suggests that the natural generalization of the Kashaev invariant
for other simple Lie algebra should be the quantum invariant of
knots colored by the Verma module of highest weight $-\delta$,
where $\delta$ is the half-sum of positive roots.

Finally we point out how the hyperbolic volume of the knot
complement, through the theory of $L^2$-torsion, has a determinant
formula that looks strikingly similar to the one of Kashaev
invariants: In both we have non-commutative deformations of the
Burau matrices, but in one case quantum determinant is use, in the
other the Fuglede-Kadison determinant is used. This  suggests an
approach to the volume conjecture using quantum determinant as an
approximation of the infinite-dimensional Fuglede-Kadison
determinant.

\subsection{A determinant formula for the colored Jones polynomial}

\subsubsection{Right-quantum matrices and quantum determinants}
A $2\times 2$ matrix $\begin{pmatrix} a & b\\
c&d
\end{pmatrix}$ is {\em right-quantum} if \begin{eqnarray*}
 ac &= & qca  \qquad  \text{(q-commutation of the entries in a column)}
 \\
 bd &= & qdb  \qquad \text{(q-commutation of the entries in a
 column)}
 \\
ad &=  & da + qcb  - q^{-1}bc \quad  \text{ (cross commutation
relation).}
\end{eqnarray*}

An $m\times m$ matrix  is {\em right-quantum} if  any $2\times 2$
submatrix of it is right-quantum. The meaning is a right-quantum
matrix preserves the structure of quantum $m$-spaces (see
\cite{Manin}). The product of 2 right-quantum matrices is a
right-quantum matrix, provided that every entry of the first
commutes with every entry of the second. The quantum determinant
of any right-quantum $A=(a_{ij})$ is defined by

$${\det}_q(A) := \sum_{\pi} (-q)^{\inv (\pi)} a_{\pi1,1} a_{\pi2,2}\dots a_{\pi m,m}$$
 where the
sum ranges over all permutations of $\{1, \dots , m\}$, and
 $\inv(\pi)$ denotes the number of inversions.

Note that in general $I-A$, where $I$ is the identity matrix, is
not right-quantum any more. We will define its determinant, using
an analog of the expansion in the case $q=1$:

$$\Det_q (I-A) := 1- C, \qquad \text{where }\quad C:=\sum_{\emptyset \neq J\subset \{1,2,\dots,r\}}
(-1)^{|J|-1}\det{}_q(A_J),$$ where $A_J$ is the $J$ by $J$
submatrix of $A$, which is always right-quantum.

\subsubsection{Deformed Burau matrix}\label{intro}
On the polynomial ring $\R[x^{\pm1},y^{\pm 1},u^{\pm1}]$ act
operators $\hat x, \tau_x$ and their inverses:
$$ \hat x f(x,y,\dots) := x f(x,y,\dots), \qquad \tau_x
f(x,y,\dots):= f(qx,y,\dots).$$

It's easy to see that $\hat x \tau_x = q \tau_x \hat x$. For other
variable, say $y$, there are similar operators $\hat y, \tau_y$,
each of which commutes with each of $\hat x, \tau_x$. Let us
define

\begin{eqnarray} a_+&=& (\hat u- \hat y\tau_x^{-1})\tau_y^{-1},
\qquad  b_+ =
\hat u^2, \qquad c_+  = \hat x \tau_y^{-2}\tau_u^{-1},  \label{t1}\\
 a_-&=& (\tau_y -\hat x^{-1})\tau_x^{-1}\tau_u,   \qquad b_- = \hat
u^{2},  \qquad c_- = \hat y^{-1}\tau_x^{-1} \tau_u. \label{t2}
\end{eqnarray}

Then it is easy to check that the following matrices $S_\pm$  are
right-quantum.

$$ S_+ := \begin{pmatrix} a & b\\
c&0
\end{pmatrix}
   \qquad  S_- :=
   \begin{pmatrix} 0 & c_-\\
b_- & a_-
\end{pmatrix}
$$

Suppose $P$ is a polynomial in the operators $a_\pm,b_\pm,c_\pm$
with coefficients in $\R=\BZ[q^{\pm1}]$. Applying $P$ to the
constant function $1$, then substituting $u$ by 1 and $x$ and $y$
by $z$, one gets a polynomial $\mathcal E(P) \in
\BZ[q^{\pm1},z^{\pm1}]$. Then it is readily seen that $\mathcal E
(S_+)$ and $\mathcal E(S_-)$ are the transpose Burau matrix and
its inverse:

$$ \mathcal E (S_+) =  \begin{pmatrix} 1-z & 1\\
z &0
\end{pmatrix}, \qquad \mathcal E (S_+) =  \begin{pmatrix} 0 & z^{-1}\\
1& 1-z^{-1}
\end{pmatrix}.$$

\subsubsection{Determinant formula} Let $\sigma_i, 1\le i \le m-1$, be the standard generators of the braid group on
$m$ strands, see for example \cite{Birman,Jones}. For a sequence
$\gamma= (\gamma_1,\gamma_2,\dots,\gamma_k)$ of pairs $\gamma_j
=(i_j,\varepsilon_j)$, where $1\le i_j \le m-1$ and
$\varepsilon_j=\pm$, let $\beta= \beta(\gamma)$ be the braid
$$\beta := \sigma_{i_1}^{\varepsilon_1}
\sigma_{i_2}^{\varepsilon_2}\dots \sigma_{i_k}^ {\varepsilon_k}.$$

Here $\sigma^\pm$ means $\sigma^{\pm 1}$. We will assume that the
closure of $\beta$ (see \cite{Birman}) is a {\em knot} , i.e. it
has only one connected component. Recall that in the Burau
representation of the braid $\beta(\gamma)$, we associate to each
$\sigma_{i_j}^{\varepsilon_j}$ an $m\times m$ matrix which is the
same as the identity matrix everywhere except for the $2\times 2$
minor of rows $i_j,i_j+1$ and columns $i_j,i_j+1$, where we put
the $2\times 2$ Burau matrix if $\ve_j=+$, or its inverse if
$\ve_j=-$. Let us do the same, only now the $2\times 2$ Burau
matrix and its inverse, for $\sigma_{i_j}^{\ve_j}$, are replaced
by $S_{+,j}$ and $S_{-,j}$. Here $S_{\pm,j}$ are the same as
$S_\pm$ with $x,y,u$ replaced by $x_j,y_j,u_j$. For the precise
definition see Section \ref{prec}. The result is a right-quantum
matrix $\rho(\gamma)$, whose entries are operators acting on
$\RR_k=\otimes_{j=1}^k\R[x_j^{\pm 1},y_j^{\pm 1},u_j^{\pm 1}]$.
Note that $\rho(\gamma)$  might not be an invariant of the braid
$\beta(\gamma)$. We can define $\mathcal E(P)$, where $P$ is an
operator acting on $\RR_k$, as before: first apply $P$ to the
constant function 1, then replace all the $u_j$ with 1, and all
the $x_j$ and $y_j$ with $z$. Further, let $\mathcal E_N(P)$ be
obtained from $\mathcal E(P)$ by the substitution $z\to q^{N-1}$.

Let $\rho'(\gamma)$ be obtained from $\rho(\gamma)$ by removing
the first row and column. Let $w(\beta)$ denotes the writhe,
$w(\beta):= \sum_j \ve_j1$. It's easy to show that when the
closure of $\beta$ is a knot, $w(\beta)-m+1$ is always even.

\begin{theorem}\label{main1}  Suppose the closure in the standard
way of the $m$-strand braid $\beta(\gamma)$ is a knot $K$.

 a) For any positive integer $N$ one has
 $$q^{(N-1)(w(\beta)-m+1)/2}\, \mathcal E_N \left (\frac{1}{\Det_q (I-q\,\rho'(\gamma))}\right)=J'_K(N).$$

b)  $\det \mathcal E(I-\rho'(\gamma))$ is equal to the Alexander
polynomial of $K$.
\end{theorem}

Part a) should be understood as follows. Suppose $\Det_q
(I-\rho'(\gamma)) = 1-C$, then when applying $\mathcal E_N$ to

\begin{equation}\frac{1}{1-C }:= \sum_{n=0}^\infty C^n,\label{d1}
\end{equation}
 only a finite number of terms are non-zero, hence the sum is
well-defined, and is equal to the colored Jones polynomial. We
would like to emphasize that here $N>0$. If $N=0$, when applying
$\mathcal E_N$ to the right hand side of (\ref{d1}), there might
be infinitely many non-zero terms. From the theorem one can
immediately get the Melvin-Morton conjecture, first proved by
Bar-Natan and Garoufalidis \cite{BNG}.
\begin{remark} Another determinant formula of the colored Jones polynomial
using non-commutative variables was given in the independent work
\cite{GLo}, also based on the quantum MacMahon Master theorem. The
main difference is here our variables are explicit operators
acting on polynomials ring. This sometimes helps since operators
can be composed. Another difference is we derive our formula from
the $R$-matrix, while \cite{GLo} used cablings of the original
Jones polynomial and graph theory. Our approach is a
non-commutative analog of Rozansky's beautiful work
\cite{Rozansky}.
\end{remark}

\subsubsection{An example} To see an application of our formula
let's calculate the colored Jones polynomial of the right-handed
trefoil. In this case we need only 2 strands with $\beta=
\sigma^3$. Thus $\rho(\gamma)= S_{+,1}S_{+,2}S_{+,3}$ is easy to
calculate, and we get $\rho'(\gamma)= c_1a_2b_3$. Hence, with $K$
being the right-handed trefoil,

\begin{eqnarray}J'_K(N) &=& q^{N-1}\,\mathcal E_N\left (\frac{1}{1-q c_1a_2b_3
}\right)=q^{N-1} \, \sum_{n=0}^\infty \mathcal E_N(q^n
c_1^na_2^nb_3^n) \notag\\
&=& q^{N-1} \sum_ {n=0}^\infty q^{n N} (1- q^{N-1})(1-q^{N-2})
\dots (1 - q^{N-n}). \label{tre}\end{eqnarray}

Note that the sum is always finite, since the term in the right
hand side is 0 if $n\ge N$.
\subsection{The Kashaev's invariant as the invariant of dimension
0} Kashaev \cite{Kashaev} used quantum dilogarithm to define a
knot invariant $\langle K\rangle_N$, depending on a positive
number $N$. Murakami and Murakami \cite{MM} showed that $\langle
K\rangle_N = J'_K(N)|_{q=\exp(2\pi i/N)}$. The famous volume
conjecture \cite{Kashaev,MM} says that the growth rate of $\langle
K\rangle_N$ is equal to the  volume $V(K)$ (see definition below)
of the knot complement:

$$\lim_{N\to \infty} \frac{\ln |\langle K\rangle_N|}{N} = \frac
{\Vol(K)}{2 \pi}.$$

 Working with varying $N$, i.e. working with varying
$sl_2$-modules might be difficult. Here we show that the values of
$\langle K\rangle_N$ comes from just one $sl_2$-module, the Verma
module of highest weight $-1$, and is a kind of analytic function
in the following sense. Let us define the
\newcommand{\Ha}{\widehat{\BZ[q]}}
Habiro ring $\Ha$ by

$$\Ha := \lim_{\leftarrow}\BZ[q]/((1-q)(1-q^2)\dots(1-q^n)).$$

Habiro \cite{Habiro} called it the cyclotomic completion of
$\BZ[q]$. Formally, $\Ha$ is the set of all series of the form

$$f(q) = \sum_{n=0}^\infty f_n(q) \, (1-q)(1-q^2)\dots(1-q^n), \qquad \text{where } \quad f_n(q) \in \BZ[q].$$

 Suppose $U$ is the set of roots of 1. If $\xi\in U$ then
 $(1-\xi)(1-\xi^2) \dots (1-\xi^n)=0$ if $n$ is big enough, hence
one can define $f(\xi)$ for $f\in \Ha$. One can consider every
$f\in \Ha$ as a function with domain $U$. Note that $f(\xi)\in
\BZ[\xi]$ is always an algebraic integer.  It turns out $\Ha$ has
remarkable properties, and plays an important role in quantum
topology. First, each $f\in \Ha$ has a natural Taylor series at
every point of $U$, and if two functions $f,g\in \Ha$ have the
same Taylor series at a point in $U$, then $f=g$. A consequence is
that $\Ha$ is an integral domain. Second, if $f=g$ at infinitely
many roots of prime power orders, then $f=g$ (see \cite{Habiro}).
Hence one can consider $\Ha$ as a class of ``analytic functions"
with domain $U$.  It was proved, by Habiro for $sl_2$ and by
Habiro with the second author for general simple Lie algebras,
that quantum invariants of integral homology 3-spheres belong to
$\Ha$ and thus have remarkable integrality properties. Here we
show that the Kashaev invariant also belongs to $\Ha$:

\begin{theorem}\label{main2}

a) $ q^{(m-w(\beta)-1)/2}\,\mathcal E_0\left (\frac{1}{\Det_q
(I-q\, \rho'(\gamma))}\right)$ belongs to $\Ha$ and is an
invariant of the knot $K$ obtained by closing $\beta(\gamma)$.

b)  Kashaev's invariant is equal to

\begin{equation}\langle K\rangle_N=  q^{(m-w(\beta)-1)/2}\,\mathcal E_0\left (\frac{1}{\Det_q (I-q\,
\rho'(\gamma))}\right)\vert_{q=\exp(2\pi i/N)}.\label{25}
\end{equation}

\end{theorem}

For example, when $K$ is the left-handed trefoil, from
\eqref{tre}, with $q\to q^{-1}$, we have

$$ \langle K\rangle_N = q\, \sum_{n=0}^\infty
(1-q)(1-q^2)\dots(1-q^n),$$ where $q=\exp(2\pi i/N)$. The function
given by the infinite sum on the right hand side was first written
down by M. Kontsevich, and its asymptotics was completely
determined by Zagier \cite{Zagier}. We see that it has a nice
geometric interpretation: It is the Kashaev invariant of the
trefoil.

\subsection{Hyperbolic volume and $L^2$-torsion}

It is known that by cutting  the knot complement $S^3\setminus K$
along some embedded tori one gets connected components which are
either Seifert-fibered or hyperbolic. Let $\Vol(K)$ be the sum of
the hyperbolic volume of the hyperbolic pieces, ignoring the
Seifert-fibered components. It's known that $\Vol(K)$ is
proportional to the Gromov norm \cite{Gromov}, and can be
calculated using $L^2$-torsion as follows. Let the knot $K$ again
be the closure of the braid $\beta$. The fundamental group  of the
knot complement has a presentation:

$$\pi_1 = \langle z_1,\dots, z_m \mid r_1,\dots,r_m\rangle,$$
where $r_i= \beta(z_i)z_i^{-1}$, with $\beta$  considered as an
automorphism of the free group on $m$ generators $z_1,\dots z_m$.

Let $Ja= \left( \frac{\partial\, r_i}{\partial\,z_j} \right)$ be
the Jacobian matrix with entries in $\BZ[\pi_1]$, where
$\frac{\partial\, r_i}{\partial\,z_j}$ is the the Fox derivative.
For a matrix with
 entries in $\BZ[\pi_1]$, one can define its Fuglede-Kadison determinant  (see \cite{Luck}), denoted by $\det_{\pi_1}$.
 A deep theorem of Luck and Schick \cite{Luck} says that

 $$\Vol(K) = 6\pi \,\ln ({\det}_{\pi_1}(Ja')),$$
where $Ja'$ is obtained from $Ja$ by removing the first row and
column. It's easy to see that

$$Ja = \psi(\beta) -I, \qquad \text{ where } \quad  \psi(\beta)=
\left( \frac{\partial(\beta(z_i))}{\partial\,z_j}\right).$$

A simple property of Fugledge-Kadison determinant is that
$\det_{\pi_1}(A) = \det_{\pi_1}(-A)$. Hence we have

\begin{proposition} Let $\psi'(\beta)$ be obtained from $\psi(\beta)$ by removing the first row and
column. Then

\begin{equation} \exp(-\frac{\Vol(K)}{6 \pi}) = \frac{1}{{\det}_{\pi_1}(I-\psi'(\beta))}.
\label{22}
\end{equation}
\end{proposition}

Note that  under the abelianization map $ab: \BZ[\pi_1] \to
\BZ[\BZ]$, the matrix $\psi(\beta)$ becomes the Burau
representation of $\beta$. Hence both $\psi(\beta)$  and
$\rho(\beta(\gamma))$  are two different kinds of quantization of
the Burau representation. We hope that the similarity between
(\ref{22}) and (\ref{25}) will help to solve the volume
conjecture. One needs to relate the Fugledge-Kadison determinant
$\det_{\pi_1}$ to the quantum determinant.

Also note that the abelianized version of the  the right hand side
of \eqref{22}, i.e. ${\det}_{\BZ}(I-ab(\psi'(\beta)))$, is equal
to the Mahler measure of the Alexander polynomial (see
\cite{Luck}). This partially explains some similarity between the
Mahler measure and the hyperbolic volume of a knot, as observed in
\cite{Silver}.

\subsection{Plan of the paper} In section \ref{proof1} we prove
Theorem \ref{main1}. Section \ref{proof2} contains a proof of
Theorem 2 and a discussion about generalization to other Lie
algebra of the Kashaev invariants.


\newcommand{\tR}{\tilde \R}

\section{Proof of Theorem \ref{main1}}\label{proof1}
In subsection \ref{pp} we recall the definition of the colored
Jones polynomial using $R$-matrix. We will follow Rozansky
\cite{Rozansky} to twist the $R$-matrix so that it has a ``nice"
form. Then in the subsequent subsections we show how the twisted
$R$-matrix can be obtained from the deformed Burau matrix, giving
a proof of Theorem \ref{main1}.

We will use the variable $v^{1/2}$ such that $v^2=q$. Note that
our $q$ is equal to $q^2$ in \cite{Jantzen}. Recall that
$\R=\BZ[q^{\pm1}]$, which is a subring of the field $\tR:=
\BC(v^{\pm 1/2})$. We will use the following standard notations.

\begin{eqnarray*}
 [n] \quad &:=& \frac{v^n-v^{-n}}{v-v^{-1}}, \qquad [n]! \quad
 :=
\prod_{i=1}^n [i], \qquad \qbinom{n}{l}\quad := \prod_{i=1}^l
\frac{[n-i+1]}{[l-i+1]},\\
(n)_{q} &:=& \frac{1-q^{-n}}{1-q^{-1}}, \qquad  \binom{n}{l}_q
\quad := \prod_{i=1}^l \frac{(n-i+1)_q}{(l-i+1)_q}, \qquad
(1-x)^d_q:= \prod_{i=0}^{d-1}(1-x q^i).
\end{eqnarray*}

\subsection{The colored Jones polynomial through $R$-matrix}
\label{pp}

\subsubsection{The quantized enveloping algebra $U_v(sl_2)$}

Let $\U$ be the algebra over the field $\tR= \BC(v^{\pm 1/2})$
generated by $K^{\pm 1/2},E,F$,  subject to the relation

$$ K^{1/2} K^{-1/2}=1,\quad K^{1/2}E= v EK^{1/2},\quad K^{1/2}F= v^{-1}FK^{1/2}, \quad EF-FE=
\frac{K-K^{-1}}{v-v^{-1}}.$$

Then $\U$ is a Hopf algebra with coproduct: $$\D(K^{1/2}) =
K^{1/2} \otimes K^{1/2}, \quad \Delta(E)= E\otimes 1 + K\otimes E,
\quad \Delta(F) = F\otimes K^{-1} + 1\otimes F.$$

 Here we follow the definition of Jantzen's book \cite{Jantzen}, only we add the square root $K^{1/2}$ for convenience.
 Note that  $V\otimes W$ has
a natural $\U$-module structure whenever $V,W$ have, due to the
co-algebra structure.

\subsubsection{The quasi-$R$-matrix and braiding}
 The quasi-$R$-matrix $\Theta$ is an
element of some completion of $\U$:

$$\Theta := \sum_{n=0}^\infty (-1)^n v^{-n(n-1)/2}\, \frac{(v-v^{-1})^n}{[n]!}\,  F^n\otimes
E^n.$$

An $\U$-module $V$ is $E$-{\em locally-finite} if for every $u\in
V$ there is $n$ such that $E^nu=0$. If $V$ and $W$ are
$E$-locally-finite, then for every $u\otimes w \in V\otimes W$,
there are only a finite number of terms in the sum of $\Theta$
that do not annihilate $u\otimes w$, hence we can define
 $\Theta$  as an $\tR$-linear operator
acting on $V\otimes W$. The inverse of $\Theta$ is given by

$$\Theta^{-1} := \sum_{n=0}^\infty  v^{n(n-1)/2}\, \frac{(v-v^{-1})^n}{[n]!}\,  F^n\otimes
E^n.$$

An element $u$ in  an $\U$-module is said to have weight $l$ if
$Ku=v^lu$. We will consider only $\U$-modules that are spanned by
weight vectors. For such modules $V$ and $W$ we define the
diagonal operator $D$ by

$$D(u\otimes w) = v^{-kl/2}u\otimes w,$$
where $u$ has weight $k$, $w$ has weight $l$. The braiding $\bb:
V\otimes W \to W\otimes V$ is defined by

$$\bb (u\otimes w) := \Theta( D (w\otimes u)).$$

It's known that $\bb$ commutes with the action of $\U$, is
invertible, and satisfies the braid relation: Suppose $V$ is an
$E$-locally-finite $\U$-module. Let $\bb_{12}:= \bb \otimes \id$
and $\bb_{23}:= \id \otimes \bb$ be the operators acting on
$V\otimes V \otimes V$. Then

$$ \bb_{12} \, \bb_{23} \, \bb_{12} = \bb_{23}
\, \bb_{12}\, \bb_{23}.$$

One can define a representation of the braid group on $m$ strands
into the group of linear operators acting on $V^{\otimes m}$ by
putting

$$ \tau(\sigma_i) = \id^{\otimes (i-1)} \otimes \,\bb \otimes
\id^{\otimes m-i-1},$$ i.e. $\sigma_i$ acts trivially on all
components, except for the $i$-th and $(i+1)$-st where it acts as
$\bb$.

\subsubsection{A modification of Verma module $V_N$}
 For an integer $N$, not necessarily positive, let $V_N$ be the $\tR$-vector space
 freely spanned by $e_i, i \in \BZ_{\ge 0}$. The following can be
 readily checked.

\begin{proposition}\label{act}  The space $V_N$ has a structure of an $E$-locally-finite $\U$-module
given by

\begin{eqnarray*}
Ke_i &= & v^{N-1-2i}e_i\\
E e_i & = & (i)_{q^{-1}}\,e_{i-1}\\
F e_i & = & v^i[N-1-i] e_{i+1}= \frac{v^{1-N}}{v-v^{-1}}\,(q^{N-1}
-q^i) \, e_{i+1}.
\end{eqnarray*}
\end{proposition}

For $N>0$ let $W_N$ be the
 $\tR$-subspace of $V_N$ spanned by $e_i, 0\le i\le N-1$. It's is easy
 to see that $W_N$ is a simple $\U$-submodule of $V_N$. Every
 simple finite dimensional $\U$-module is isomorphic to one of
 $W_N$.

 \begin{remark}The traditional basis $e'_i := F^i(e_0)/[i]!$ is
 related to the basis $e_i$ by

 $$ \qbinom{N-1}{i}\, e_i = v^{-i(i-1)/2}\, e'_i.$$

 \end{remark}

\subsubsection{The colored Jones polynomial}
If the closure of the $m$-strand braid $\beta$ is the knot $K$,
then the colored Jones polynomial $J_K(N)$ can be defined as the
quantum trace of $\tau(\beta)$ on $(W_N)^{\otimes m}$:

$$ J_K(N) = v^{w(\beta) \frac{N^2-1}{2}}\tr_q(\tau(\beta), (W_N)^{\otimes m})
:= v^{w(\beta) \frac{N^2-1}{2}}\tr(\tau(\beta)K^{-1}, (W_N)^{\otimes m}).$$

Here $w(\beta) := \sum_j \ve_j 1$ is the writhe of $\beta$. The
factor $v^{w(\beta) \frac{N^2-1}{2}}$ will make $J_K(N)$ not
depending on the framing. If $K$ is the unknot then $J_K(N)= [N]$.
The normalized version $J'_K(N) := J_K(N)/[N]$ can be calculated
using the partial trace
 as follows. Recall that $\tau(\beta)$ acts on $(W_N)^{\otimes
 m}$. Taking the quantum trace of $\tau(\beta)$ in only the $m-1$ last components, we get an
 operator acting on the first $W_N$, which is known to be a scalar times the
 identity operator, with  the scalar being exactly $J'_K(N)$. This can be
 written in the formula form as follows.
Let  $p_0: (V_N)^{\otimes m} \to (V_N)^{\otimes m}$ be the
projection onto $e_0\otimes (V_N)^{\otimes (m-1)}$, i.e.
$$p_0(e_{n_1}\otimes e_{n_2}\otimes \dots \otimes e_{n_m})=
\delta_{0,i_1}\, e_{n_2}\otimes \dots \otimes e_{n_m}.$$

Then $p_0$ also restricts to a projection from $(W_N)^{\otimes m}$
onto $e_0\otimes (W_N)^{\otimes (m-1)}$, and

\newcommand{\ptr}{\operatorname {ptr}}

\begin{equation} J'_K(N) = v^{w(\beta) \frac{N^2-1}{2}}
 \tr\left( p_0(\tau(\beta)\, K^{-1}) , e_0 \otimes (W_N)^{\otimes (m-1)}\right).
\label{ptr}
\end{equation}

\subsubsection{Twisting the braiding} It's straightforward to calculate the
action of the braiding $\bb$ on $V_N\otimes V_N$, using the basis
$e_{n_1}\otimes e_{n_2}, n_1,n_2\in \BZ_{\ge 0}$. However to get a
better, more convenient form we will follow Rozansky
\cite{Rozansky} to use the twisted braiding

$$ \check \bb  := Q^{-1} \bb \,Q \qquad  \text{ where }\quad Q= \id\otimes
K^{(1-N)/2}.$$

Then direct calculation shows that on $V_{N}\otimes V_{N}$ the
action of the twisted braiding $\check \bb _\pm$ are given by

$$ \check \bb _\pm (e_{n_1}\otimes e_{n_2}) =  \sum_{l=0}^{\max n_1,
n_2} \check \bb _\pm (n_1,n_2;l) \, (e_{n_2\pm l} \otimes
e_{n_1\mp l}),$$ where, with $z=q^{N-1}$,

\begin{eqnarray}
(\check \bb _+)(n_1,n_2;l)& = & q^{-\frac{(N-1)^2}{4}}\,
\binom{n_1}{l}_{q^{-1}}\, q^{n_2(l-n_1)} z^{n_2} \, (1-z q^{-n_2})^l_{q^{-1}} \label{a8}\\
(\check \bb _-)(n_1,n_2;l)& = & q^{\frac{(N-1)^2}{4}}\,
\binom{n_2}{l}_{q}\, q^{n_1(n_2-l)} z^{-n_1}\,  (1-z^{-1}
q^{n_1})^l_{q}. \label{a9}
 \end{eqnarray}

Note that our formulas differ from those in \cite{Rozansky} by
$q\to q^{-1}$, since we derived our formula directly from the
quantized enveloping algebra that differs from the one implicitly
used by Rozansky. (The co-products are opposite; ``implicitly"
since Rozansky never used quantized enveloping algebra, but just
took the formula of the $R$-matrix from \cite{KM}).

To justify the use of the twisted braiding we argue as follows.
First note that $\bb_\pm $ commutes with $K^{1/2}$, the action of
which on $V_N\otimes V_N$ is given by $\D(K^{1/2}) =
K^{1/2}\otimes K^{1/2}$. Thus $K^l\otimes K^l$ commutes with
$\bb_\pm $ for every half-integer $l$. Hence

\begin{equation}(Q')^{-1} \, \bb_\pm \, Q' = Q^{-1} \, \bb_\pm \, Q= \check
\bb _\pm \label{oi} \end{equation}
 if $$Q' = Q \, (K^{i(1-N)/2}\otimes
K^{i(1-N)/2})= K^{i(1-N)/2}\otimes K^{(i+1)(1-N)/2}.$$

Let us define the operator $Q_m$ acting on $(W_N)^{\otimes m}$ by
$$Q_m:= K^{(1-N)/2} \otimes K^{2(1-N)/2} \otimes \dots\otimes
K^{m(1-N)/2}$$ and let
$$\check \tau(\beta) =Q_m^{-1} \tau(\beta)
Q_m.$$

 Then $\check \tau$ is also a representation of the braid
group. Since the action of $K^{-1}$ on $(W_N)^{\otimes m}$
commutes with the action of $Q_m$, one sees that in the formula
(\ref{ptr}) we can use $\check \tau(\beta)$ instead of
$\tau(\beta)$:

\begin{equation}
J'_K(N) = v^{w(\beta) \frac{N^2-1}{2}} \tr\left( p_0(\check
\tau(\beta)\, K^{-1}) , e_0 \otimes (W_N)^{\otimes (m-1)}\right).
\label{ptr2}
\end{equation}

Suppose  $\beta = \sigma_{i_1}^{\ve_1} \dots
\sigma_{i_k}^{\ve_k}$. Then $\check \tau(\beta) = \check
\tau(\sigma_{i_1})^{\ve_1} \dots \check
\tau(\sigma_{i_k})^{\ve_k}$. Let us calculate $\check
\tau(\sigma_{i})$:

\begin{eqnarray*}
\check \tau(\sigma_{i}^{\pm 1} ) &=&  Q_m^{-1} \, \tau(\sigma_i)\,  Q_m\\
&=& Q_m^{-1} \,(\id^{\otimes (i-1)} \otimes \,\bb_\pm \otimes
\id^{\otimes m-i-1}) \, Q_m\\
&=& \id^{\otimes (i-1)} \otimes \, \left( (K^{i(1-N)/2}\otimes
K^{(i+1)(1-N)/2})^{-1} \, \bb_\pm \, (K^{i(1-N)/2}\otimes
K^{(i+1)(1-N)/2}) \right) \otimes \id^{\otimes m-i-1}\\
&=& \id^{\otimes (i-1)} \otimes \, \check \bb _\pm \otimes
\id^{\otimes m-i-1} \qquad \text{ by  \eqref{oi}}.
\end{eqnarray*}

This means in the definition of $\check \tau$ one just use $\check
\bb _\pm$ instead of $\bb_\pm$, and then $\check\tau$ is obtained
from $\tau$ by the global twist $Q_m$.

\subsubsection{From $W_N$ to $V_N$} So far we take the trace using the finite dimensional module $W_N$.
For the infinite dimensional $V_N$ we define the trace of an
operator if only a finite number of diagonal entries are nonzero.
The following was observed in \cite{Rozansky}.

\begin{lemma} Suppose the closure of the braid $\beta$ is a knot,
then
\begin{eqnarray*}
J'_K(N) &=& v^{w(\beta) \frac{N^2-1}{2}}\, \tr\left( p_0(\check
\tau(\beta)\, K^{-1}) , e_0 \otimes (W_N)^{\otimes (m-1)}\right)
\\
&= & v^{w(\beta) \frac{N^2-1}{2}} \,\tr\left( p_0(\check
\tau(\beta)\, K^{-1}) , e_0 \otimes (V_N)^{\otimes (m-1)}\right).
\end{eqnarray*}
\label{a1}
\end{lemma}

\begin{proof}

One important observation is that if $n < N$, and $n+l\ge N$, then
$F^l e_n=0$. Hence $\check \bb _\pm (e_{n_1}\otimes e_{n_2})$ is a
linear combination of $e_{m_1}\otimes e_{m_2}$ (with
$m_1+m_2=n_1+n_2$), and if $n_1 <N$ then $m_2 <N$, or if $n_2 <N$
then $m_1 < N$.

Let $\left (\tau(\beta)K^{-1}\right)_{n_1,n_2\dots
n_m}^{s_1,s_2\dots s_m}$ be the matrix of $\tau(\beta)K^{-1}$ with
respect to the basis $e_{n_1}\otimes e_{n_2}\otimes \dots \otimes
e_{n_m}$ in $(V_N)^{\otimes m}$. Note that $K^{-1}$ acts
diagonally in this basis. The above observation shows that if $n_i
<N$ then $s_{\bar \beta(i)} <N$ for the matrix entry $\left
(\tau(\beta)K^{-1}\right)_{n_1,n_2\dots n_m}^{s_1,s_2\dots s_m}$
not to be 0, where $\bar \beta$ is the permutation corresponding
to $\beta$. To take the trace we only have to concern with the
case $s_i=n_i$. We have already had $n_1=0$, which is less than
$N$. Thus we must have $n_j <N$ for $j=1, \bar\beta(1),
(\bar\beta)^2(1)\dots$. The fact  that the closure of $\beta$ is a
knot implies that $\{(\bar\beta)^l(1), 1\le l \le m\}$ is the
whole set $\{ 1, 2, \dots, m\}$. Hence taking the trace over $e_0
\otimes (V_N)^{\otimes (m-1)}$ is the same as over $e_0 \otimes
(W_N)^{\otimes (m-1)}$.
\end{proof}

\subsection{Algebra of the deformed Burau matrix}\label{ppp}

\subsubsection{Algebra $\A_\ve$}
Let us define

\begin{eqnarray*}
\A_+ & := \R\langle a_+, b_+, c_+\rangle /(a_+b_+ = b_+ a_+,
a_+c_+
 = q c_+a_+,  b_+c_+ = q^2 c_+b_+).\\
\A_- &:= \R\langle a_-, b_-, c_-\rangle /(a_-b_- = q^{2} b_- a_-,
c_-a_-  = q a_-c_-,   c_-b_- = q^2b_-c_-).
\end{eqnarray*}

It is easy to check that the $a_\pm,b_\pm,c_\pm$ of section
\ref{intro} satisfy the commutation relations of the algebras
$\A_\pm$.

 For a sequence $\varepsilon= (\varepsilon_1,\varepsilon_2,\dots,\varepsilon_k)$,
where each $\varepsilon_j$ is either $+$ or $-$,  let
$\A_\varepsilon = \A_{\varepsilon_1}\otimes
\A_{\varepsilon_2}\otimes \dots \otimes \A_{\varepsilon_k}$. We
can consider $\A_\varepsilon$ as the algebra over $\R$ freely
generated by $a_j,b_j,c_j$ subsect to the commutation relations:
if $i\neq j$ then each of $a_i,b_i,c_i$ commutes with each of
$a_j,b_j,c_j$, if $\varepsilon_j=+$ then the commutations among
$a_j,b_j,c_j$ are the same as those of $a_+,b_+,c_+$, and if
$\varepsilon_j=-$ then the commutations among $a_j,b_j,c_j$ are
the same as those of $a_-,b_-,c_-$.  Note that the algebra
$\A_\varepsilon$ is a generalized quantum space in the sense that
for any $a,b$ among the generators, one has the almost
$q$-commutation relation $ab = q^{l} ba$, for some integer $l$.

Replacing $x,y,u,a_\pm,b_\pm,c_\pm$  with respectively
$x_j,y_j,u_j, a_j,b_j,c_j$ in (\ref{t1}) if $\varepsilon_j=+$, or
in (\ref{t2}) if $\varepsilon_j=-$, we identify $a_j,b_j,c_j$ with
operators acting on $\R[x_j^{\pm1}, y_j^{\pm1},u_j^{\pm1}]$. We
assume that $a_j,b_j,c_j$  leave alone $x_i,y_i,u_i$ if $i\neq j$.
Thus $\A_\varepsilon$ acts on the algebra $\P_k$ of Laurent
polynomials in $x_j,y_j,u_j, 1\le j\le k$ with coefficients in
$\R$. The map $\mathcal E : \A_\varepsilon \to \R[z^{\pm1}]$ is
defined as in section \ref{intro}.

\begin{lemma} a) If $f,g\in \A_\ve$ are separate, i.e.  $f$ contains only
$a_j,b_j,c_j$ with $j \le r$ and $g$ contains  only $a_l,b_l,c_l$
with $r<l$ (for some $r$), then $\mathcal E(fg) = \mathcal E(f)
\mathcal E(g)$.

b) One has \begin{eqnarray}
\mathcal E(b_+^s c_+^r a_+^d) & = &q^{-rd} \, z^r \, (1-z q^{-r})^d_{q^{-1}} \label{s1}\\
\mathcal E(b_-^s c_-^r a_-^d) & = &z^{-r} \, (1-z^{-1}
q^{r})^d_{q}\label{s2}
\end{eqnarray}

\label{adic}
\end{lemma}
\begin{proof} a) follows directly from the definition. b) follows
from an easy induction.\end{proof}

\subsubsection{Definition of $\rho(\gamma)$}\label{prec}
 Let us give here the precise definition of
$\rho(\gamma)$, for $\gamma=((i_1,\ve_1),\dots,(i_k,\ve_k))$.
Recall that $\beta$ is the braid
$$ \beta= \beta(\gamma):=  \sigma_{i_1}^{\varepsilon_1}
\sigma_{i_2}^{\varepsilon_2}\dots \sigma_{i_k}^ {\varepsilon_k}.$$

If $\ve_j=+$ (resp. $\ve_j=-$), let  $S_{j}$ be the matrix $S_+$
(resp. $S_-$) with $a_+,b_+,c_+$ (resp. $a_-,b_-,c_-$) replaced by
$a_{j},b_{j},c_{j}$. For the $j$-th factor
$\sigma_{i_j}^{\varepsilon_j}$ let us define an $m\times m$
right-quantum matrix $A_j$ by the block sum, just like in the
  Burau representation, only the non-trivial $2\times 2$
block now is $ S_{j}$ instead of the Burau matrix:

$$ A_j  := I_{i_j-1} \oplus S_{j}\oplus
I_{m-i_j-1}.$$ Here $I_l$ is the identity $l\times l$ matrix.

Let $\rho(\gamma) := A_1 A_2\dots A_k$. Then $\rho(\gamma)$ is an
$m\times m$ right-quantum matrix with entries polynomials in
$a_{j}, b_{j},c_{j}$.

\subsection{Quantum MacMahon Master Theorem}

\subsubsection{Co-actions of right-quantum matrices on the quantum space}
The quantum plane $\BC_q[z_1,z_2, \dots, z_m]$, considered as the
space of $q$-polynomial in the variables $z_1,\dots, z_m$, is
defined as
$$\BC_q[z_1,z_2, \dots, z_m]:= \tR\langle
z_1,\dots,z_m\rangle/(z_iz_j=qz_jz_i \text{ if } \, i<j).$$

\begin{remark}Our definitions of quantum spaces, quantum matrices... differ from the one in \cite{GLZ,Kassel} by
the involution $q\to q^{-1}$, but agree with the ones in Jantzen's
book \cite{Jantzen}.
\end{remark}

 If $A=(a_{ij})_{i,j=1}^m$ is right-quantum and all $a_{ij}$'s commute
with all $z_1,\dots,z_m$, then it is known that the $Z_i:= \sum_j
a_{ij} z_j$, i.e.
$$\begin{pmatrix} Z_1\\ Z_2\\ \dots \\ Z_m\end{pmatrix} = A \begin{pmatrix} z_1\\ z_2\\ \dots \\ z_m\end{pmatrix},$$
also satisfy $Z_iZ_j=qZ_jZ_i \text{ if } \, i<j$.
\newcommand{\AAs}{\mathcal W}
Let $\AAs=\AAs(A)$ be the algebra generated by $a_{ij}, 1\le i,j
\le m$, subject to the commutation relations of $a_{ij}$.  Then we
have an {\em algebra} homomorphism:

$$\Phi_A: \BC_q[z_1\dots,z_m] \to \AAs \otimes \BC_q[z_1,\dots,z_m]$$
defined by $\Phi_A(z_i)=Z_i$. Informally, one could look at
$\Phi_A$ as the degree-preserving {\em algebra} homomorphism on
the $q$-polynomial ring $\BC_q[z_1,z_2, \dots, z_m]$ defined by
matrix $A$. Here we assume that the degree of each $z_i$ is 1, and
the degree of each $a_{ij}$ is 0.

We will consider the case $A=\rho(\gamma)$, and in particular $A=
\check \bb _\pm$. In this case we define $\mathcal E_N(\Phi_A):=
(\mathcal E_N\otimes id)\, \circ \, \Phi_A$, which is a {\em
linear} operator acting on $\BC_q[z_1\dots,z_m]$, not necessarily
an algebra homomorphism.

\subsubsection{Quantum MacMahon Master theorem} Let
$\BC_q[z_1,\dots,z_m]^{(n)}$ be the part of total degree $n$ in
$\BC_q[z_1,\dots,z_m]$. Since $\Phi_A$ preserves the total degree,
it restricts to a linear map: $ \Phi_A: \BC_q[z_1\dots,z_m]^{(n)}
\to \AAs \otimes \BC_q[z_1,\dots,z_m]^{(n)}$. Let us define the
trace by

$$\tr\left(\Phi_A, \BC_q[z_1\dots,z_m]^{(n)}\right)=
\sum_{n_1+\dots+n_m=n}(\Phi_A)_{n_1,\dots,n_m}^{n_1,\dots,n_m},$$
where $(\Phi_A)_{n_1,\dots,n_m}^{n_1,\dots,n_m}$ is the
coefficients of $z_1^{n_1}\dots z_m^{n_m}$ in $Z_1^{n_1}\dots
Z_m^{n_m}$. One could consider $\tr\left(\Phi_A,
\BC_q[z_1\dots,z_m]^{(n)}\right)$ as the trace of $\Phi_A$ acting
on the part of total degree $n$. The quantum MacMahon's Master
theorem, proved in \cite{GLZ} says that

$$\frac{1}{\Det_q (I-A)} = \tr\left(\Phi_A,
\BC_q[z_1\dots,z_m]\right):= \sum_{n=0}^\infty \tr\left(\Phi_A,
\BC_q[z_1\dots,z_m]^{(n)}\right).$$

It's the $q$-analog of the identity

$$\frac{1}{\det (I-C)}= \sum_{n=0}^\infty \tr(S^nC),$$
where $C$ is a linear  operator acting on a finite dimensional
$\BC$-space $V$ and $S^nC$ is the action of $C$ on the $n$-th
symmetric power of $V$.

\subsection{From deformed Burau matrices $S_\pm$ to $R$-matrices $\check \bb _\pm$}

Let $\F_m:(V_N)^{\otimes m} \to \BC_q[z_1,\dots,z_m]$ be the
$\tR$-linear isomorphism defined by $\F(e_{n_1}\otimes \dots
\otimes e_{n_m}) := z_1^{n_1}\dots z_m^{n_m}$. The following is
important to us.

\begin{proposition}\label{a2}

a) Under the isomorphism $\F_2$, the twisted braiding matrices
$\check \bb _\pm$ acting on $V_N\otimes V_N$ map to  $ v^{\mp
(N-1)^2/2}\, \mathcal E_N(S_\pm)$, i.e.
$$ \check \bb _\pm =  v^{\mp
(N-1)^2/2}\, \F_2^{-1}\, \mathcal E_N(\Phi_{S_\pm}) \, \F_2.$$

b) Under the isomorphism $\F_m$, the linear automorphism $\check
\tau(\beta(\gamma))$ of $(V_N)^{\otimes m}$ maps to $ v^{\mp
w(\beta) (N-1)^2/2}\mathcal E_N(\Phi_{\rho(\gamma)})$.

\end{proposition}

\begin{proof}

a) Suppose for 2 variables $X,Y$ we have $YX=qXY$, then Gauss's
q-binomial formula \cite{Kassel} says that

$$ (X+Y)^n = \sum_{l=0}^n \binom{n}{l}_q X^l Y^{n-l}.$$

Let us first consider the case of $S_+$. Then $\Phi_{S_+}(z_1)=
a_+ z_1 + b_+ z_2$, and $\Phi_{S_-}(z_2) = c_+ z_1$. Note that
$(b_+ z_2) (a_+z_1)=q^{-1}(a_+z_1) b_+(z_2)$, hence using the
Gauss binomial formula we have

\begin{eqnarray*} \Phi_{S_+}(z_1^{n_1} z_2^{n_2} )&=& (a_+
z_1 + b_+z_2)^{n_1} (c_+z_1)^{n_2}\\
&=& \sum_{l=0}^{n_1}  \binom{n_1}{l}_{q^{-1}} (a_+z_1)^l\, (b_+
z_2)^{n_1-l}
\, (c_+z_1)^{n_2}\\
&=& \sum_{l=0}^{n_1}  \binom{n_1}{l}_{q^{-1}} \, q^{n_2(n_1-l)}\,
a_+^l\, b_+^{n_1-l} c_+^{n_2} \,  (z_1)^{n_2+l}\, (z_2)^{n_1-l}
\end{eqnarray*}

Using formulas (\ref{a8}) and (\ref{s1}) one sees that

$$ \check \bb _+ =  v^{- (N-1)^2/2}\, \F_2^{-1}\, \mathcal
E_N(\Phi_{S_+}) \, \F_2.$$

The proof for $S_-$ is quite similar, using formulas (\ref{a9})
and (\ref{s2}).

b) Because the variables $x_j,y_j,u_j$ are separated, we have that

$$\mathcal E(\rho(\gamma)) =
\mathcal E(\rho(\sigma_{i_1}^{\ve_1})\dots \mathcal
E(\rho(\sigma_{i_k}^{\ve_k}),$$ and the statement follows from
part a).
\end{proof}

\subsubsection{} Under the isomorphism $\F_m$, the projection
$p_0: (V_N)^{\otimes m}\to (V_N)^{\otimes m}$ maps to the
projection, also denoted by $p_0$, of $\BC_q[z_1, z_2,\dots,z_m]$,
which can be defined as

$$p_0(z_1^{n_1}\, z_2^{n_2}\dots z_m^{n_m})= \delta_{0,n_1}\, z_2^{n_2}\dots
z_m^{n_m}.$$

Note that the kernel of $p_0$ is the ideal generated by $z_1$.

\begin{lemma} a) For every $u\in \BC_q[z_2,\dots,z_m]$,
$$p_0  \left (\Phi_{\rho(\beta)}(u)\right) = \Phi_{\rho'(\beta)}(u).$$

b) The operators $p_0$ and $\mathcal E_N$ commute:

$$ p_0  \left( \mathcal E_N\left (\Phi_{\rho(\beta)}(u)\right)\right)=
\mathcal E_N\left (p_0  \left
(\Phi_{\rho(\beta)}(u)\right)\right).$$ \label{a3}
\end{lemma}
\begin{proof}
a) Recall that $\rho'(\gamma)$ is obtained from $\rho(\gamma)$ by
removing the first row and column. Suppose $u \in
\BC_q[z_2,\dots,z_m]$, then $\Phi_{\rho(\beta)}(u) -
\Phi_{\rho'(\beta)}(u)$ is divisible by $z_1$, and hence
annihilated by $p_0$.

b) follows trivially from the definition. \end{proof}

The following is trivial.

\begin{lemma} Under $\F_m$, the action of $K^{-1}$ on
$\BC_q[ z_2,\dots,z_m]^{(n)}$ is the scalar operator, with scalar
$v^{(m-1)(1-N)+ 2n} = v^{(m-1)(1-N)} q^n$. \label{a6}
\end{lemma}

\subsection{Proof of Theorem \ref{main1}}

\begin{eqnarray*}
J'_K(N) &=& v^{\frac{w(\beta)(N-1)^2}{2}}\,\tr\left( p_0
\left(\check\tau(\beta) K^{-1}\right), e_0\otimes
(W_N)^{\otimes (m-1)} \right) \qquad \text{by Lemma \ref{a1}}\\
&=&  v^{w(\beta)(N-1)}\, \tr\left( p_0 \left( \mathcal
E_N(\Phi_{\rho(\gamma)}) K^{-1}\right) ,
\BC_q[z_2,\dots,z_m] \right) \qquad \text{under $\F_m$, by Proposition \ref{a2}}\\
&= & v^{w(\beta)(N-1)}\, \tr \left( \mathcal E_N(\Phi_{\rho'(\gamma)})K^{-1}, \BC_q[z_2,\dots,z_m] \right)   \quad \text{by Lemma \ref{a3}}\\
&=& v^{w(\beta)(N-1)}\, \sum_{n=0}^\infty \tr\left( \mathcal
E_N(\rho'(\gamma))
K^{-1}, \BC_q[z_2,\dots,z_m]^{(n)} \right)\\
&=& v^{(w(\beta)-m+1)(N-1)}\, \sum_{n=0}^\infty q^n \, \tr\left(
\mathcal E_N(\Phi_{\rho'(\gamma)}),
\BC_q[z_2,\dots,z_m]^{(n)} \right)  \quad \text{by Lemma \ref{a6}} \\
&=& v^{(w(\beta)-m+1)(N-1)}\, \sum_{n=0}^\infty \tr\left( \mathcal
E_N( \Phi_{q\rho'(\gamma)}),
\BC_q[z_2,\dots,z_m]^{(n)} \right)\\
&=& v^{(w(\beta)-m+1)(N-1)}\, \mathcal E_N \sum_{n=0}^\infty
\tr\left( \Phi_{q\rho'(\gamma)},
\BC_q[z_2,\dots,z_m]^{(n)} \right)\\
&=& v^{(w(\beta)-m+1)(N-1)}\,  \mathcal E_N  \frac{1}{\Det_q(I-q\,
\rho'(\gamma))} \qquad \text{by quantum MacMahon Master Theorem.}
\end{eqnarray*}

\newcommand{\leb}{\overset{\leftarrow}{\beta}}

This proves part a) of Theorem \ref{main1}. As for part b), first
notice that the braid $\leb:= \sigma_{i_k}^{\ve_k}\sigma
_{i_{k-1}}^{\ve_{k-1}}\dots \sigma_{i_1}^{\ve_1}$ has the closure
knot the same as that of $\beta$. The Alexander polynomial of $K$
is known to be equal to $\det(I-\bar \rho'(\leb))$, where
$\bar\rho$ is the Burau representation, and $\bar \rho'(\leb)$ is
obtained from $\bar \rho(\leb)$ by removing the first row and
column. We know that $\mathcal E(S_\pm)$ are the transpose Burau
matrices, hence $\bar \rho(\leb)= \mathcal E(\rho(\beta))^T$, the
transpose of $\mathcal E(\rho(\beta))$. The statement now follows.

\section{The Kashaev invariant}

\subsection{Proof of Theorem \ref{main2}}\label{proof2}
 \subsubsection{Completion of $\A_\ve$}

Let $\I$ be the left ideal in $\A_\varepsilon$ generated by
$a_1,a_2,\dots,a_k$, i.e.
$$ \I := a_1\A_\varepsilon +
a_2\A_\varepsilon + \dots + a_k \A_\varepsilon,$$
 and let $\hat
\A_\varepsilon$ be the $\I$-adic completion of $\A_\varepsilon$.
Using the almost $q$-commutativity it's easy to see that $\I$ is a
two-sided ideal.

\begin{lemma} When the closure of $\beta(\gamma)$ is a knot,   $\Det_q (I-q\rho'(\gamma))$ belongs to $1+\I$, and
hence $\frac{1}{\Det_q(I-q\rho'(\gamma))}$ belongs to $\hat
\A_\ve$. \label{adic2}
\end{lemma}

\begin{proof} It's enough to show that when $a_1=a_2=\dots a_k=0$,
then $\Det_q(I-q\rho'(\gamma))=1$, or $\det_q(C)=0$ for any main
minor $C$ of $\rho'(\gamma)$.

Let call {\em permutation-like} matrix a matrix $C$ where on each
row and on each column there is at most one non-zero entry. If, in
addition, on each row and on each column there is exactly one
non-zero entry, we say that $C$ is non-degenerate. Every
non-degenerate permutation-like square matrix $C$ gives rise to a
permutation matrix $p(C)$ by replacing all the non-zero entries
with 1. It's clear that product of (non-degenerate)
permutation-like matrices is a (non-degenerate) permutation-like
one. If $C$ is a permutation-like $m\times m$ matrix, and $D$ a
main minor, i.e. a submatrix  of type $J\times J$, then $D$ is
also permutation like. If, in addition, both $C$ and $D$ are
non-degenerate, then $p(C)$ leaves $J$ stable, i.e. $p(C)(J)=J$,
since the restriction of $p(C)$ on $J$ is equal to $p(D)$ which
leaves $J$ stable.

Also note that if $C$ is degenerate permutation-like right-quantum
matrix, then $\det_q(C)=0$.

When $a_1=a_2=\dots =a_k=0$, each of matrices $A_j$ (whose
definition is in subsection \ref{prec}) is a non-degenerate
permutation-like matrix. Hence $C=\rho(\gamma)$ is
permutation-like. Note that $p(C)$ is exactly $\bar\beta$, the
permutation corresponding to $\beta$. Because the closure of
$\beta$ is a {\em knot}, $\bar\beta=p(C)$ does not leave any
proper subset of $\{1,2,\dots,m\}$ stable. Hence any main minor
$D$ of $\rho'(\gamma)$, which itself is a proper main minor of $C=
\rho(\gamma)$, is a degenerate permutation-like matrix. Hence
$\det_q(D)=0$.\end{proof}

\subsubsection{$\hat \A_\ve$ and the Habiro ring}

\begin{lemma} a) If $f\in \A_\ve$ is divisible by $a_j^d$ for some $1\le j\le k$ and a positive integer $d$, then
$\mathcal E(f)$ is divisible by $(1-zq^{r})^d_q$, and hence
$\mathcal E_N(f)$ is divisible by $(1-q)^d_q$ for every integer
$N$, not necessarily positive.

b) Suppose $n >dk$. Then   $\mathcal E_N(f)$is divisible by
$(1-q)^d_q$ for every integer $N$ and every $f\in \I^n$. Hence
$\mathcal E_N \hat A_\varepsilon \in \Ha$. \label{adic1}
\end{lemma}
\begin{proof} a)
We assume that $f$ is a {\em monomial} in the variables
$a_1,b_1,c_1,a_2,\dots$. Using the almost $q$-commutativity we
move all $a_j,b_j,c_j$ to  the right of $f$, so that $f=g \, b_j^s
c_j^r a_j^d$, for some $g\in \A_\ve$ not containing $a_j,b_j,c_j$.
Note that by Lemma \ref{adic}
$$\mathcal E(f) =\mathcal E(g)\,
\mathcal E(b_j^s c_j^r a_j^d) $$ is divisible by $\mathcal E(b_j^s
c_j^r a_j^d)$. Note that $a_j,b_j,c_j$ are either $a_+,b_+,c_+$ or
$a_-,b_-,c_-$. Using \eqref{s1} and \eqref{s2} we see that
$\mathcal E_N(f)$ is divisible by $(1-q^l)^d_q$ for some integer
$l$, which, in turn, is always divisible by $(1-q)^d_q$.

 b) Using the fact that generators $a_j,b_j,c_j, 1\le j\le k$ almost $q$-commute, it's easy to see that
 $\I^n$ is 2-sided ideal generated by $a_{s_1}a_{s_2}\dots a_{s_n}$, where each $s_i$ is one of $\{1,2,\dots,k\}$.
 If $n>dk$, by the pigeon hole principle, there is an index $j$ such $a_{s_1}a_{s_2}\dots
 a_{s_n}$ is divisible by $a_j^d$. Now the result follows from
 part a). \end{proof}

 From Lemmas \ref{adic1} and \ref{adic2} we get the following.

\begin{corollary} Suppose $N$ is an integer, not necessarily
positive. Then

$$ \mathcal E_N \left(\frac{1}{\Det_q(I-q\rho'(\gamma))}\right) \in
\Ha.$$ \label{cor}
\end{corollary}

\subsubsection{Proof of Theorem \ref{main2}}

Part a) is a special case of Corollary \ref{cor}, with $N=0$.

 For part b) first recall that $\langle
K\rangle _N = J'_K(N)|_{q=\exp(2\pi i/N)}$. When $q=\exp(2\pi
i/N)$, one has $q^N=1=q^0$. Thus $\mathcal E_N=\mathcal E_0$ when
$q=\exp(2\pi i/N)$. One has

\begin{eqnarray*}J'_K(N)|_{q=\exp(2\pi i/N)} &=&  v^{m-1-w(\beta)}\, \mathcal
E_N (T)|_{q=\exp(2\pi i/N)}\\
&= &  v^{m-1-w(\beta)}\, \mathcal E_0(T)|_{q=\exp(2\pi
i/N)},\end{eqnarray*} where
$$T=\frac{1}{\Det_q(I-q\rho'(\gamma))}.$$

\subsection{The Kashaev invariant for other simple Lie algebra}

Fix a simple Lie algebra $\g$. For every long knot K, presented by
a $1-1$ tangle, one can define the $\g$-universal invariant
$\J_K^{\g}$, which is a central element in an appropriate
completion of quantized universal enveloping algebra $U_v(\g)$,
see \cite{Turaev,Lawrence}. Formally, $\J_{K,\g}$ is an infinite
sum of central elements in $U_v(\g)$:

\begin{equation}
\J_{K,\g} = \sum_{n=0}^\infty \J_{K,\g}^{(n)}, \label{o1}
\end{equation}
such that for any finite dimensional simple $U_v(\g)$-module only
the action of a finite number of terms are non-zero. Hence for a
finite-dimensional simple module $U_v(\g)$-module $V$, $\J_K^{\g}$
acts as a scalar times the identity. It can be shown that the
scalar is a Laurent polynomial in $q$. Denote this scalar by
$J'_{K,\g}(V)$. One always has

$$J_{K,\g}(V)= J'_{K,\g}(V)\dim_q(V),$$
 when $J_{K,\g}(V)$ is the
usual quantum invariant of $K$ colored by $V$, and $\dim_q(V)$ is
the quantum dimension, i.e. the invariant of the unknot colored by
$V$.

For any Verma module $V_\lambda$ of highest weight $\lambda$ (an
element in the weight lattice), the action of each of
$\J_{K,\g}^{(n)}$ is still in $\R=\BZ[q^{\pm1}]$, but in general
infinitely many of them are non-zero. In this case
$J'_{K,\g}(V_\lambda)$ is an infinite series (sum). In a future
work we will show that $J'_{K,\g}(V_\lambda)\in \Ha$; the special
case when $\g=sl_2$ has been proved here by
 Corollary \ref{cor}.

Note that if the weight $\lambda$ is dominant, then

$$ J'_{K,\g}(V_\lambda)= J'_{K,\fg}(W_\lambda),$$ where $W_\lambda$
is the finite dimensional $U_v(\fg)$ module with highest weight
$\lambda$. The reason is both are the scalar of the same scalar
operator acting on $V_\lambda$ and its quotient $W_\lambda$. In
this case $ J'_{K,\g}(V_\lambda)$  is a Laurent polynomial in $q$.
It is known that $\R=\BZ[q^{\pm 1}]\subset \Ha$, see
\cite{Habiro}.

Due to the Weyl symmetry, we see that if $w$ is in the Weyl group,
then $J'_{K,\g}(V_\lambda) = J'_{K,\g}(V_{w\cdot \lambda})$, where
$w\cdot \lambda$ is the dot action of the Weyl group, see
\cite{Humphreys}. If $\lambda$ is not fixed (under the dot action)
by any element of the Weyl group, then $\lambda = w\cdot \mu$ for
some dominant $\mu$, and hence  $J'_{K,\g}(V_\lambda)=
J'_{K,\g}(V_\mu)$. In this case $J'_{K,\g}(V_\lambda)$ might be
still an infinite series, but it is equal to a Laurent polynomial,
which is $J'_{K,\fg}(V_\mu)$ in the Habiro ring $\Ha$.

The more interesting, and less understood case is when $\lambda$
is fixed by an element of the Weyl group, i.e. $\lambda$ is on a
wall of a {\em shifted} Weyl chamber. Among them there is one
special weight, namely $\lambda= -\delta$, where $\delta$ is the
half-sum of positive roots, since $-\delta$ is the only element
invariant by {\em all} elements of the Weyl group. When
$\fg=sl_2$, $V_{-\delta}$ is $V_0$ in section \ref{proof1}, and
$J'_{K,\g}(V_{-\delta})$ is the Kashaev invariant in this case,
according to Theorem \ref{main2}. Note that $V_{-\delta}$ is
always infinite-dimensional and irreducible; it's certainly a very
special $U_v(\fg)$-module.

Thus a natural generalization of the Kashaev invariant to other
simple Lie algebra is $J'_{K,\g}(V_{-\delta})$. More precisely,
let's define the $\fg$-Kashaev invariant by

$$ \langle
K\rangle _N^{\g} := J'_{K,\g}(V_{-\delta})|_{q=\exp(2\pi i/N)}.$$

And we suggest the following $\fg$-volume conjecture

$$\lim_{N\to \infty} \frac{|\langle
K\rangle _N^{\g}|}{N} = c_\fg \, \Vol(K),$$

Where $c_\fg$ is a constant depending only on the simple Lie
algebra $\fg$.

\end{document}